\documentclass{elsart}
\usepackage{amssymb}
\usepackage{amsmath}
\usepackage[english]{babel}

\newtheorem{theorem}{Theorem}

\newtheorem{lemma}[theorem]{Lemma}

\newtheorem{remark}[theorem]{Remark}

\begin{document}
\begin{frontmatter}
\title{Asymptotic approximation of degenerate fiber integrals.}
\author{Brice Camus}
\thanks[label2]{This work was supported by the
\textit{SFB/TR12} project, \textit{Symmetries and Universality in
Mesoscopic Systems}.}

\address{Ruhr-Universit\"at Bochum, Fakult\"at f\"ur Mathematik,\newline%
Universit\"atsstr. 150, D-44780 Bochum, Germany.\newline%
Email : brice.camus@univ-reims.fr}

\begin{abstract}
We study asymptotics of fiber integrals depending
on a large parameter. When the critical fiber is singular, full-asymptotic expansions are established in two different
cases : local extremum and isolated real principal type singularities. The main coefficients are computed
and invariantly expressed. In the most singular cases it is
shown that the leading term of the expansion is related to
invariant measures on the spherical blow-up of the singularity. The results can be applied
to certain degenerate oscillatory integrals which occur in spectral analysis and quantum mechanics.
\end{abstract}
\end{frontmatter}
\section{Introduction and statement of the main result.}
In \cite{B-S} J.Br\"uning\& R.Seeley have studied asymptotic
expansions of integrals :
\begin{equation}
H(z)=\int\limits_{0}^{\infty} \sigma (xz,x)dx,\text{ }
z\rightarrow \infty,
\end{equation}
where $\sigma(x,\xi)$ is a singular symbol. The result of \cite{B-S} is quite remarkable, in particular because it
can be directly applied to spectral analysis. Many asymptotic questions can be reduced to the study of the
previous problem but it is also interesting to consider a
generalization :
\begin{equation}\label{the problem}
I(z)=\int\limits_{X} g(z f(x),x)dx,\text{ }
z\rightarrow \infty.
\end{equation}
where $g: \mathbb{R}\times X \rightarrow \mathbb{R}$, $f:X\rightarrow \mathbb{R}$
are smooth and $X$ is a smooth differentiable manifold equipped with the $C^{\infty}$ density $dx$.
As a motivation, we remark that the trace formula for
certain linear self-adjoint operators can be stated in the form of Eq. (\ref{the
problem}). In the semi-classical regime when $f$ is a quadratic form this problem was involved in \cite{BPU}.
See also \cite{Cam1} for the case of $f$ with an isolated and degenerate
singularity associated to an homogeneous definite jet.

\textit{General assumptions.} Throughout this work, we will
assume that $|f|$ is strictly positive outside of a compact set
and that the Fourier transform $\hat{g}$ w.r.t. $t$ exists with
$\partial_t^k \hat{g}(t,x)\in L^1(\mathbb{R}\times X)$, $\forall k$.
For this reason $g$ will be called a symbol. $\hfill{\square}$
\begin{remark}\rm{This assumption on $g$ is strong but can be weakened.
More general conditions can be found in \cite{HOR1}. Mainly, this condition will be used
to reach integrals with compact supports.}
\end{remark}
As $z\rightarrow \infty$, the asymptotic behavior of $I(z)$ is related to the critical fiber :
\begin{equation}
\mathfrak{S}=f^{-1}(\{0\})=\{ x\in X \text{ / } f(x)=0\}.
\end{equation}
This can easily be viewed with the Fourier inversion formula :
\begin{equation}\label{problem osci}
I(z)= \int\limits_{X} \int\limits_{\mathbb{R}} e^{i z t f(x)}
\hat{g}(t,x)dtdx,\text{ } z\rightarrow \infty,
\end{equation}
where $\hat{g}(t,x)$ is the normalized Fourier transform of $g$ w.r.t. $t$ :
\begin{equation}
\hat{g}(t,x)=\frac{1}{2\pi}\int\limits_{\mathbb{R}} e^{-i\tau t}
g(\tau,x)d\tau,
\end{equation}
In Eq. (\ref{problem osci}), the stationary points w.r.t. $t$ are
precisely given by $\mathfrak{S}$, i.e. $I(z)$ is asymptotically
supported by $\mathfrak{S}$. Since $f$ is smooth $\mathfrak{S}$ is closed and
according to the general assumptions above we obtain :
\begin{equation}
(H_0) \textit { The fiber }\mathfrak{S} \textit{ is compact}.
\end{equation}
$(H_0)$ is not absolutely necessary, e.g. if $g$
decreases fast enough near the boundary of $X\cap \mathfrak{S}$. But, to simplify, we will
only consider the compact case. To obtain an easier formulation of the problem we recall
an elementary result.
\begin{lemma}\label{non stationary}
If $\partial^k_t\hat{g}\in L^1(\mathbb{R}\times X)$, $\forall k\in\mathbb{N}$, modulo terms $\mathcal{O}(z^{-\infty})$,
asymptotics of Eq. (\ref{the problem})
are not changed by assuming that
$\hat{g}$ is compactly supported near
$\mathfrak{S}$.
\end{lemma}\label{compact-form}
\textit{Proof.} With $\mathfrak{S}$ compact we choose a cut-off
function $\Theta\in C_0^{\infty}(X)$ such that $\Psi=1$ near
$\mathfrak{S}$ and $0\leq \Psi \leq1$. We shall estimate
the error integral :
\begin{equation}
E(z)=\int\limits_{X} \int\limits_{\mathbb{R}} e^{i z t f(x)}
\hat{g}(t,x)(1-\Psi(x))dtdx,\text{ } z\rightarrow \infty.
\end{equation}
With $L=-(i/zf(x))\partial_t$, we have
$L^k e^{i z t f(x)}=e^{i z t f(x)}$, $\forall k \in\mathbb{N}$.
By integration by parts and since $|f(x)|\geq C$ on $\mathrm{supp}(1-\Psi)$, we obtain :
\begin{equation}
|E(z)|\leq (Cz)^{-k} || \partial_t^{k} \hat{g}(t,x)
||_{\mathrm{L}^1(\mathbb{R}\times X)}=\mathcal{O} (z^{-k}), \text{ } \forall k
\in\mathbb{N}.
\end{equation}
This gives the desired result, with our hypothesis on $g$.
$\hfill{\blacksquare}$\medskip\\
Lemma \ref{compact-form} allows to consider only integrals with compact support
w.r.t. $x$ which simplifies all questions of convergence. Notice that we can weaken the condition
on $g$ to $\partial^k_t\hat{g}\in L^1(\mathbb{R}\times X)$, $\forall k\leq k_0$ with an error $\mathcal{O}(z^{-k_0})$.
We are mainly interested in the situation where $\mathfrak{S}$ has an isolated singularity.
If $x_0\in\mathfrak{S}$ is such a critical point, let $\Theta\in C_0^\infty(X)$ be a cut-off microlocally
supported near $x_0$. We split-up our integral as
$I(z)=I_r(z)+I_s(z)$, where :
\begin{gather}
I_r(z)=\int\limits_{X}  g(z f(x),x)(1-\Theta)(x)dx, \label{regularpart}\\
I_s(z)=\int\limits_{X} g(z f(x)t,x) \Theta(x)dx. \label{singpart}
\end{gather}
This procedure can be extended with finitely many critical points on $\mathfrak{S}$.
The regular part $I_r$ can be treated by the generalized stationary phase method,
with non-degenerate normal Hessian, which we recall in section 2.
Since $I_s$ is a local object and the main
contributions below concern invariant objects, there is no
loss of generality to assume that $\mathrm{supp}(\Theta)$ is an open of
$\mathbb{R}^n$, $n=\dim (X)$. For $x_0\in \mathfrak{S}$ a singularity of finite order we can write the germ of $f$ as :
\begin{equation}\label{germ}
f(x)=f_k(x)+\mathcal{O}(||(x-x_0)||^{k+1}),
\end{equation}
where $f_k\neq 0$ is homogeneous of degree $k\geq 2$ w.r.t.
$(x-x_0)$. The first elementary result concerns extremum attached to such
homogeneous germs :
\begin{theorem} \label{main1}
If $f$ has a local extremum $x_0$ on $\mathfrak{S}$ whose jet is given by Eq.
(\ref{germ}) (a fortiori $k$ is even), we obtain a full-asymptotic expansion :
\begin{equation}
I_s(z)\sim\sum\limits_{j\in\mathbb{N}} c_j z^{-\frac{j}{k}}.
\end{equation}
If $\dim(X)=n$, the leading term is given by :
\begin{equation} \label{sing part for extremum}
I_s(z)=z^{-\frac{n}{k}} \left\langle t_e^{\frac{n-k}{k}}, g(t,x_0)
\right\rangle \frac{1}{k}\int\limits_{\mathbb{S}^{n-1}}
|f_k(\theta)|^{-\frac{n}{k}} d\theta
+\mathcal{O}(z^{-\frac{n+1}{k}}),
\end{equation}
with $t_e=\max(t,0)$ if $x_0$ is a minimum and $\max(-t,0)$ for a
maximum.
\end{theorem}
The reader can observe that Theorem \ref{main1} includes the case of a non-integrable
singularity on $\mathfrak{S}$ for $k>n$. Accordingly, always for $k>n$,
the contribution of the critical point is bigger than the regular contribution (see section 2).
\begin{remark} \rm{The result stated in Eq. (\ref{sing part for extremum}) includes the case
$k=2$ and the integral on the sphere can be expressed in terms of
special functions. Under ad-hoc conditions, Theorem \ref{main1} can be extended to a singularity which is a sum
of positively homogeneous singularities, e.g. $f(x_1,x_2)=||x_1||^4+||x_2||^6$.
}
\end{remark}
The case of non-extremum degenerate critical points is more difficult. Since this problem can be
very complicated in general position we impose :\medskip\\
\textit{$(H_1)$ $\mathfrak{S}$ has a unique critical point $x_0$. Moreover, $f_k$ defined in Eq. (\ref{germ}) is
non-degenerate in the sense that :}
\begin{equation}\nabla f_k \neq 0 \text{ }on \text{ }\mathbb{S}^{n-1}\cap
\{f_k=0\}.
\end{equation}
By homogeneity, $(H_1)$ is equivalent to "$f_k$ has an isolated
singularity". Observe that this condition is very close to H\"ormander's real principal type condition
for distributions. We define the integrated density of $f_k$ on the sphere (or co-area) as :
\begin{equation}
\mathrm{LVol}(w)=\int\limits_{\{f_k(\theta)=w\}} |dL|(\theta),\text{ }dL(\theta)\wedge df_k(\theta)=d\theta.
\end{equation}
Where $dL$ is the $n-2$ dimensional Liouville measure induced by $f_k$ on $\mathbb{S}^{n-1}$,
i.e. the Riemannian density induced by $f_k$ on the standard density of $\mathbb{S}^{n-1}$.
Note that $(H_1)$ insures that $\mathrm{LVol}(w)$ is well defined and smooth near the origin.
\begin{theorem}\label{main2}
Under the previous assumptions and if $x_0$ satisfies $(H_1)$, the
singular part of our integral admits a full asymptotic expansion :
\begin{equation}
I_s(z)\sim \sum\limits_{j=0}^{\infty} c_{j}z^{-\frac{j}{k}}
+\sum\limits_{j=0}^{\infty} d_j z^{-j}\log(z).
\end{equation}
a) If $k>n$ (non-integrable singularity), the leading term is :
\begin{equation}
I_s(z)=C_0 z^{-\frac{n}{k}}
+\mathcal{O}(z^{-\frac{n+1}{k}}\log(z)),
\end{equation}
where the distributional coefficient $C_0$ is given by :
\begin{equation*}
\frac{1}{k}(\left\langle t_{+}^{\frac{n}{k}-1}, g(t,x_0)\right\rangle
\int\limits_{\{f_k \geq 0\}}
|f_k(\theta)|^{-\frac{n}{k}}d\theta  +\left\langle
t_{-}^{\frac{n}{k}-1}, g(t,x_0)\right\rangle
\int\limits_{\{f_k \leq 0\}}
|f_k(\theta)|^{-\frac{n}{k}}d\theta).
\end{equation*}
b) If $k$ divides $n$ the leading term is logarithmic :
\begin{equation}
I_s(z)= D_0 z^{-\frac{n}{k}} \log (z)
+\mathcal{O}(z^{-\frac{n}{k}}),
\end{equation}
with $n=kp$, we obtain :
\begin{equation*}
D_0=\frac{1}{k}
\left ( \frac{d^{p-1}}{dw^{p-1}}\mathrm{Lvol}(0) \right ) \int\limits_{\mathbb{R}} |t|^{p-1}g(t,x_0)dt.
\end{equation*}
c) If $n>k$ and $n/k \notin \mathbb{N}$ (integrable singularity) we obtain the same result as in a) but
with the modified distributions :
\begin{equation*}
\left\langle t_{+}^{\frac{n}{k}-1}, g(t,x_0)\right\rangle
\left\langle \frac{\tilde{d}^n}{\tilde{d}w^n} w^{n-\frac{n}{k}}_{+}, \mathrm{Lvol} \right\rangle
+\left\langle t_{-}^{\frac{n}{k}-1}, g(t,x_0)\right\rangle
\left\langle \frac{\tilde{d}^n}{\tilde{d}w^n} w^{n-\frac{n}{k}}_{-}, \mathrm{Lvol} \right\rangle,
\end{equation*}
where the derivatives w.r.t. $w$ are normalized distributional derivatives.
\end{theorem}
The meaning of normalized derivative it that one choose the normalization :
\begin{equation}\label{normalized derivatives}
\left\langle \frac{\tilde{d}^n}{\tilde{d}w^n} w^{n-\frac{n}{k}}_{\pm}, f(w)\right\rangle :=
\left\langle w^{-\frac{n}{k}}_{\pm}, f(w)\right\rangle,
\end{equation}
for all $f\in C_0^\infty$ with $f=0$ in a neighborhood of the origin.
The distributional bracket involving $\mathrm{Lvol}$ is detailed in section 3.
Results c) and b) for $p\geq 2$ are not intuitive and are certainly difficult to be
reached without geometry. In particular, for applications to oscillatory integrals (see below) one has to work
in the dual since both Fourier transforms w.r.t. $t$ in c) and b) are distributional.
In c), the $n$-th derivative is arbitrary and the result is the same for any normalized derivative of order greater than
$\mathrm{E}(n/k)$.
\begin{remark} \rm{Results a) and b) for $p=1$ are certainly interesting for spectral analysis
since these contributions are bigger than $I_r(z)=\mathcal{O}(z^{-1})$. As in Theorem \ref{main1},
non-integrable singularities have a dominant contribution. Hence, the leading term of $I(z)$ is always an invariant.}
\end{remark}
\textbf{Application to oscillatory integrals.}\\
A typical application of Theorems \ref{main1}\&\ref{main2} can be the asymptotic expansion of
distributional traces of quantum propagators. Hence, it is interesting to remark that our results
can be extended to asymptotic integrals :
\begin{equation*}
\tilde{I}(z)=\int\limits_{X} G(z,zf(x),x)dx,\text{ }z\rightarrow
+\infty,
\end{equation*}
if $G$ admits an asymptotic expansion with a priori estimates,
i.e. :
\begin{gather*}
G(z,t,x)=\sum\limits_{j=0}^{l} z^{-\alpha_j}g_j(t,x)+R_l(z,t,x),\\
\forall k\in\mathbb{N}^* : ||R_k(z,t,x)||_{L^1(\mathbb{R}\times
X)}=\mathcal{O}(z^{-(\alpha_k+\varepsilon)}),\text{ } \varepsilon >0,
\end{gather*}
where $(\alpha_j)_j$ is a strictly increasing sequence and $R$
controlled by uniform estimates. Similarly, we can consider
expansions in term of $z^{-\alpha_j}\log(z)^m$. This notion of graduation w.r.t. $z$ allows to
apply our results but, to simplify the exposition, in this work
we will just consider the case of an integral of a symbol
$g(t,x)$. We can treat degenerate oscillatory integrals :
\begin{equation}\label{IO}
O(z)=\int\limits_{\mathbb{R}\times X} e^{iz t f(x)} a(t,x,z)dtdx, \text{ }z\rightarrow +\infty,
\end{equation}
providing that $f$ satisfies the conditions of Theorem \ref{main1} or \ref{main2}.
An important application in quantum mechanics is the case $X=T^*\mathbb{R}^n$ where, after some technical
modifications, the localized (distributional) trace of $h$-pseudors :
\begin{equation*}
\mathrm{Tr } u_h (A_h-E):= \mathrm{Tr} \int\limits_{\mathbb{R}}  \hat{u}(t) e^{\frac{i}{h}t(A_h-E)} dt,
\text{}  \hat{u}\in C_0^\infty(\mathbb{R}), \text{ } E\in\mathbb{R},
\end{equation*}
can be written as a locally finite sum of oscillatory integrals :
\begin{equation}\label{IO quantum}
\int\limits_{\mathbb{R}\times T^*\mathbb{R}^n} e^{\frac{i}{h}( S(t,y,\eta)-\left\langle y,\eta \right\rangle -tE)}
b(h,t,y,\eta) dtdyd\eta,
\end{equation}
where $b(h,\bullet)\sim \sum h^{-k} b_k$ satisfies a priori estimates as above and $S$ is the local generating
function of the group of diffeomorphism of the principal symbol of $A_h$.
Here $z=h^{-1}$ is the parameter and, after a
discussion based on classical mechanics, Eq. (\ref{IO quantum}) can be reformulated as in Eq. (\ref{IO}) where
$\mathfrak{S}$ is the energy surface of level $E$. For more details, we refer to \cite{BPU,Cam0,Cam1}.

We recall now basics on homogeneous transformations, some of them will also be used below.
The Melin transform of a function $h$ is defined as :
\begin{equation}
M[h](\xi)=\int\limits_{0}^{\infty} h(t) t^{\xi-1}dt.
\end{equation}
Generally $M[h](\xi)$ is analytic, for
example, in the strip $\Re(\xi)\in]a,b[$ with :
\begin{equation*}
a= \inf\limits_{x\in\mathbb{R}} \{ \int\limits_{0}^\infty
|h(t)|t^{x-1} dt <\infty \},\text{ }
b=\sup\limits_{x\in\mathbb{R}} \{ \int\limits_{0}^\infty
|h(t)|t^{x-1} dt <\infty \},
\end{equation*}
and when $M[h](c+iy)\in L^1(\mathbb{R},dy)$ we have the inversion
formula :
\begin{equation}
h(t)=\frac{1}{2i\pi} \int\limits_{c-i\infty}^{c+i\infty} t^{-\xi}
M[h](\xi)d\xi.
\end{equation}
By elementary changes of path we obtain that :
\begin{equation*}
M[e^{\pm i t}](\xi)= e^{\pm i \pi \frac{\xi}{2}}\Gamma(\xi),\text{ }
\mathrm{Re}(\xi)\in ]0,1[.
\end{equation*}
We write  $O(z)=O_{+}(z)+O_{-}(z)$, where :
\begin{equation*}
O_{\pm}(z)=\int\limits_{\{\pm tf(x)>0\}} e^{iztf(x)}a(t,x)dtdx.\\
\end{equation*}
Melin's inversion formula leads to the distributional formulation :
\begin{equation*}
O_{\pm}(z)=\frac{1}{2i\pi} \int\limits_{c-i\infty}^{c+i\infty}
e^{\pm i\xi\frac{\pi}{2}} \Gamma(\xi) z^{-\xi}
\int\limits_{\{\pm tf(x)>0\}}
|t f(x)|^{-\xi} a(t,x)dtdx d\xi,
\end{equation*}
where $0<c<\delta$ for $\delta>0$ depending on $f$.
If no elementary method works, one can construct explicit meromorphic extensions
of distributions :
\begin{gather*}
\left\langle T_{f}^{+}(\xi),a\right\rangle=e^{+i\xi\frac{\pi}{2}} \Gamma(\xi) z^{-\xi}\int\limits_{\{tf(x)>0\}}
(tf(x))^{-\xi} a(t,x)dxdt,\\
\left\langle T_{f}^{-}(\xi),a\right\rangle=e^{-i\xi\frac{\pi}{2}} \Gamma(\xi) z^{-\xi}\int\limits_{\{tf(x)<0\}}
|tf(x)|^{-\xi} a(t,x)dxdt,
\end{gather*}
and use Cauchy's residue formula to get asymptotic w.r.t.
$z\rightarrow \infty$. This method is described in \cite{WON}
when $f$ is monomial, which is a generic situation when $f$ is
analytic on $\mathrm{supp}(a)$ by Hironaka's theorem of resolution
of singularities \cite{Hir}. Also in many cases it is possible to reach this
setting by a non-analytic change of coordinates. This approach is very interesting but a great disadvantage
is that there is 4 domains, a sequence of poles (of order 3) which does not contribute at all in the expansion and
the complex factors $\exp(\pm i\pi \xi /2) \Gamma(\xi)$ lead to long calculations in presence of multiple poles.\\
In fact we will reach Theorem \ref{main1}\&\ref{main2} more directly and by a method which avoids
Fourier analysis until the last step of the proof.
\section{Case of $\mathfrak{S}$ regular and non-degenerate critical points.}
If $f$ is regular on $\mathfrak{S}$, then $\mathfrak{S}$ is a compact and smooth submanifold.
\begin{prop}Under the previous assumptions and if $\mathfrak{S}$ is a regular
surface, $I(z)$ admits a full asymptotic expansion in powers of
$z^{-1}$ with :
\begin{equation}
I(z)=\frac{1}{z} \int\limits_{x\in\mathfrak{S}}
\int\limits_{\mathbb{R}}g(t,x)
dtd_{\mathfrak{S}}(x)+\mathcal{O}(z^{-2}),
\end{equation}
where $d_{\mathfrak{S}}$ is the invariant surface measure of
$\mathfrak{S}$. The same result holds for the integral $I_r(z)$
with insertion of the cut-off in the integral.
\end{prop}
By the implicit functions theorem and compactness, we pick some
coordinates $y$, defined near $\mathfrak{S}$, and a finite
partition of unity $\Omega_j$ covering $\mathfrak{S}$, such that
$f$ is diffeomorphic to the coordinate $y_1$ in
$\mathrm{supp}(\Omega_j)$. By Lemma \ref{non stationary} we obtain
:
\begin{gather*}
I(z)=\sum I_j(z)+\mathcal{O}(z^{-\infty}),\\
I_j(z)=\int\limits_{\mathbb{R}^2}\int\limits_{\mathbb{R}^{n-1}}
e^{i z t y_1} y^{*}(\Omega_j(x) \hat{g}(t,x)|Jy|) dtdy_1dy_2
...dy_n,
\end{gather*}
where $y^{*}$ and $|Jy|$ are respectively the pullback by $y$ and
the multiplication by the standard Jacobian. Hence, we have to estimate an oscillatory integral with
quadratic phase. The stationary phase Lemma, see e.g. Lemma
7.7.3 of \cite{HOR1} vol.1, provides a full-asymptotic expansion :
\begin{equation}
\int\limits_{\mathbb{R}^2}e^{i z t y_1} y^{*}(\Omega_j(x)
\hat{g}(t,x)|Jy|) dtdy_1=\sum\limits_{l=0}^{N-1} C_l
z^{-(l+1)} +\mathcal{O}(z^{-(N+1)}).
\end{equation}
Since $\mathrm{sign}(ty_1)=0$, the expansion is real with leading term :
\begin{equation}
I_j(z)=\frac{2\pi}{z} \int\limits_{\mathbb{R}^{n-1}}
(y^{*}(\Omega_j(x) \hat{g}(t,x)|Jy|))(0,0,y_2,...,y_n)
dy_2...dy_n+\mathcal{O}(z^{-2}).
\end{equation}
Geometrically, since $\mathfrak{S}$ is locally given by $y_1=0$,
this corresponds to integration on $ \mathfrak{S}$ w.r.t. the
canonical $n-1$ dimensional measure $d_{\mathfrak{S}}$. By
summation over the indices $j$ we obtain the integral on the whole
surface, i.e. :
\begin{equation}\label{density-regul.}
I(z)\sim \frac{2\pi}{z} \int\limits_{x\in\mathfrak{S}} \hat{g}(0,x)
d_{\mathfrak{S}}(x)=\frac{1}{z} \int\limits_{x\in\mathfrak{S}}
\int\limits_{\mathbb{R}}g(t,x)
dtd_{\mathfrak{S}}(x).
\end{equation}
\begin{remark}\rm{Eq. (\ref{density-regul.}) is a particular (flat) case
of stationary phase formulas with a smooth compact manifold of critical points
and non-degenerate transverse Hessian. $d_{\mathfrak{S}}$ is the Liouville-measure of classical
mechanics or Guelfand-Leray-measure in theory of singularities.
The oscillatory representation of delta-Dirac
distributions, by mean of Schwartz kernels, provides a very natural definition of this object.}
\end{remark}
Now, assume that $f$ admits a single non-degenerate critical point
on $\mathfrak{S}$. As usually, non-degenerate means that
$Q(\xi)=\frac{1}{2}\left\langle \xi,d^2f(x_0)\xi\right\rangle$ is a non-degenerate quadratic form. After perhaps a
reduction of the cut off $\Theta$, the Morse Lemma, with coordinates $w$, provides a simpler problem :
\begin{equation} \label{IO cubic}
F(z)=\int\limits_{\mathbb{R}\times\mathbb{R}^n} e^{izt Q(w)}
w^{*} (\hat{g}(t,x) \Theta(x)|Jw|(x)) dw.
\end{equation}
This kind of asymptotic problem was precisely studied in
\cite{BPU}. The result is :
\begin{prop}
If $Q$ has at least one even index of inertia then :
\begin{equation}
F(z)\sim \sum\limits_{j=0}^{\infty} c_j z^{-1-\frac{j}{2}}.
\end{equation}
But, if $Q$ has both indices of inertia odd then we have :
\begin{equation}
F(z) \sim \sum\limits_{j=0}^{\infty} \sum\limits_{m=0,1} c_{j,m}
z^{-1-j}\mathrm{log}(z)^m.
\end{equation}
\end{prop}
The method used in section 3.2 allows to find this result. Our
proof is not simpler and that's why we refer to Proposition 3.4
and Theorem 3.5 of \cite{BPU} for a detailed proof and discussion
of the coefficients. This discussion is necessary because the exponential distributions
of Eq. (\ref{IO cubic}) possess spherical symmetries.
\section{Proof of the main results.}
\subsection{Critical points attached to a local extremum.}
We start by the case where $f$ admits a local minimum $x_0$ on the
fiber $\mathfrak{S}$. We split-up our integral $I(z)$ via Eqs. (\ref{regularpart},\ref{singpart}).
$I_r$ can be treated as in section 2 and from now we concentrate our attention on $I_s$.
With the extremum condition, $x_0$ is isolated on $\mathfrak{S}$
and $\mathfrak{S}\cap\mathrm{supp}(\Theta)=\{x_0\}$ for $\mathrm{supp}(\Theta)$ small
enough.To simplify notations we identify $x_0$ with the origin and we recall that we consider :
\begin{equation}
f(x)=f_{k}(x)+\mathcal{O}(||x||^{k+1}),
\end{equation}
where $f_k$ is homogeneous of degree $k$ and definite positive.
Using polar coordinates, the Taylor formula shows that :
\begin{equation}
f(r\theta)=r^k (f_k(\theta)+R(r,\theta)),
\end{equation}
with $R(0,\theta)=0$ and where we note again $f_k(\theta)$ the
restriction of $f_k$ to $\mathbb{S}^{n-1}$. If $\mathrm{supp}(\Theta)\subset B(0,r_0)$
is chosen small enough we have :
\begin{equation*}
f_k(\theta)+R(r,\theta)\neq 0,\text{ } \forall
\theta\in\mathbb{S}^{n-1}, \text{ } \forall r \in [0,r_0[.
\end{equation*}
Hence, with the homogenous coordinates $v=(u,\theta)$ where :
\begin{equation}
u(r,\theta)=r (f_{k}(\theta)+R(r,\theta))^{\frac{1}{k}},
\end{equation}
we can express our integral as :
\begin{equation}\label{reduced form mini}
I_s(z)=\int\limits_{\mathbb{R}_{+}} G(z u^k,u)du.
\end{equation}
The new symbol $G$ is obtained by pullback and integration :
\begin{equation}\label{amplitude 1}
G(t,u)=\int\limits_{\mathbb{S}^{n-1}} v^*(g(t,r\theta)\Theta(r\theta)r^{n-1}|J(v)|)d\theta.
\end{equation}
The existence of the asymptotic expansion is a consequence of :
\begin{lemma}\label{lemma for extremum}
For $a$ in
$C_0^{\infty}(\mathbb{R}\times\mathbb{R}_{+})$, the following asymptotic expansion holds :
\begin{equation}
J(z)=\int\limits_{\mathbb{R}_{+}}
a(zu^k,u)du \sim \sum\limits_{j\geq 0}
z^{-\frac{j+1}{k}}d_j(a),\text{ }z\rightarrow \infty,
\end{equation}
where the $d_j$ are universal distributions given by :
\begin{equation}
d_j=\frac{1}{k}
\frac{1}{j!}(\tau_{+}^{\frac{j+1-k}{k}}\otimes
\delta_0^{(j)}(u)).
\end{equation}
\end{lemma}
\textit{Proof.} The crucial point is that $u^k$ is increasing on $\mathbb{R}_+$. We have :
\begin{equation*}
J(z)=z^{-\frac{1}{k}}
\int\limits_{0}^{\infty}
a(u^k,\frac{u}{z^{\frac{1}{k}}})du.
\end{equation*}
A Taylor expansion w.r.t. the second argument at the origins gives
:
\begin{equation}
a(u^k,\frac{u}{z^{\frac{1}{k}}})=\sum\limits_{l=0}^{N}
\frac{1}{l!} z^{-\frac{l}{k}} u^l \frac{\partial^l
a}{\partial u^l}(u^k,0)+z^{-\frac{N+1}{k}}R_{N+1}(u,z),
\end{equation}
where $R_{N+1}(u,z)$ is integrable w.r.t. $u$ with $L^1$ norm
uniformly bounded in $z$. By a new change of variable we obtain :
\begin{equation}
J(z)=\frac{1}{k} \sum\limits_{l=0}^{N} \frac{1}{l!}
z^{-\frac{1+l}{k}} \int\limits_{0}^{\infty} \frac{\partial^l
a}{\partial u^l}(\tau,0) \tau^{\frac{l+1-k}{k}} d\tau
+\mathcal{O}(z^{-\frac{N+1}{k}}).
\end{equation}
These coefficients are well defined since
$|\tau|^{\frac{l+1-k}{k}}\in L^1_{\mathrm{loc}}(\mathbb{R})$ for all
$l\in\mathbb{N}$. $\hfill{\blacksquare}$
\begin{remark}\label{remark maximum} \rm{This expansion holds also for a pullback by
$-u^k$ if we replace $\tau_{+}^{\alpha}$ by $\tau_{-}^{\alpha}$.
This allows to treat the case of a local maximum in Theorem
\ref{main1}.} Also for an application to oscillatory integrals we obtain a nice formulation
via the Fourier transform of the distributions $\tau_{\pm}^{\alpha}$, which avoids any "regularization".
\end{remark}
\textit{Proof of Theorem \ref{main1}.} We treat the case of a
local minimum. We apply Lemma \ref{lemma for extremum} to the
integral of Eq. (\ref{reduced form mini}) to prove the existence
of the asymptotic expansion and it remains to express
invariantly the leading term. With the polar coordinates $G$ vanishes up to the
order $n-1$. Consequently, we have :
\begin{equation*}
I_s(z)=\frac{z^{-\frac{n}{k}}}{k} \frac{1}{(n-1)!} \left\langle
t_{+}^{\frac{n-k}{k}}\otimes \delta_0 ^{n-1},
G \right\rangle+\mathcal{O}(z^{-\frac{n+1}{k}}).
\end{equation*}
Starting form Eq. (\ref{amplitude 1}) and since :
\begin{equation*}
|Jv|(0,\theta)=|f_k(\theta)|^{-\frac{1}{k}},
\end{equation*}
by elementary manipulations on delta-Dirac distributions, we obtain that :
\begin{equation}\label{formula for C}
I_s(z)= z^{-\frac{n}{k}}\left\langle t_{+}^{\frac{n-k}{k}},g(t,0)\right\rangle
\frac{1}{k}\int\limits_{\mathbb{S}^{n-1}}|f_k(\theta)|^{-\frac{n}{k}}
d\theta+\mathcal{O}(z^{-\frac{n+1}{k}}).
\end{equation}
Finally, for a local maximum we replace the distributions $t_{+}^{\frac{n-k}{k}}$ by $t_{-}^{\frac{n-k}{k}}$.
$\hfill{\blacksquare}$
\begin{remark}
\rm{The leading term is invariantly defined
since Eq. (\ref{formula for C}) involves the evaluation $g(t,z_0)$
and $f_k$ is an invariant since the $(k-1)$-jet in $z_0$ is flat.}
\end{remark}
\textbf{On the integral on the sphere.}\medskip\\
The integrals on $\mathbb{S}^{n-1}$ of Eq. (\ref{formula for C}) can be
reformulated. For example in dimension 2 these are elliptic
integrals. In higher dimensions we define :
\begin{equation}
I(f_k)=\int\limits_{\mathbb{R}^n} e^{-|f_k|(x)}dx,
\end{equation}
since $f_k$ is homogeneous we obtain :
\begin{equation*}
I(f_k)=\int\limits_{\mathbb{R}_{+}\times \mathbb{S}^{n-1}} e^{-r^k
f_k(\theta)} r^{n-1} dr d\theta=\int\limits_{0}^{\infty} e^{-u^k}
u^{n-1} du \int\limits_{\mathbb{S}^{n-1}}
|f_k(\theta)|^{-\frac{n}{k}} d\theta.
\end{equation*}
Hence our integral is given by :
\begin{equation}
\frac{1}{k}\int\limits_{\mathbb{S}^{n-1}}
|f_k(\theta)|^{-\frac{n}{k}} d\theta=\frac{I(f_k)}{\Gamma (n/k)}.
\end{equation}
$I(f_k)$ can be computed, by elementary methods, as a product of
gamma factors or hypergeometric functions. We give an elementary
example in section 4.
\subsection{Case of non-extremum critical points.}
This problem is more complicated and we recall that we just treat
here singularities given by $(H_1)$. As before, we just have to
study $I_s(z)$ and, in local coordinates, we identify
the critical point with the origin. We define $C(f_k)$ as the trace on the unit sphere of the conical set
of the zeros of $f_k$, i.e. :
\begin{equation}
C(f_k)=\{ \theta \in \mathbb{S}^{n-1} \text{ / } f_k(\theta)=0\}
=\mathbb{S}^{n-1}\cap \{x\in\mathbb{R}^n \text{ / } f_k(x)=0\}.
\end{equation}
With $(H_1)$, $C(f_k)$ is a compact smooth submanifold of $\mathbb{S}^{n-1}$ of codimension 1.
To perform a blow-up of the singularity we use polar
coordinates $(r,\theta)$ and the next lemma gives a resolution of the
singularity w.r.t. $C(f_k)$.
\begin{lemma}\label{normal form}
In a micro-local neighborhood of the origin there exists local
coordinates $y$, on the blow-up of the critical point, such that :
\begin{gather*}
f(x) \simeq \left\{
\begin{matrix}
y_{1}^{k},\text{ in all directions where }f_{k}(\theta )>0,\\
-y_{1}^{k},\text{ in all directions where }f_{k}(\theta )<0,\\
y_{1}^{k}y_{2},\text{ locally near }C(f_{k}).
\end{matrix}
\right.
\end{gather*}
\end{lemma}
\noindent\textit{Proof.} By Taylor, there exists $R$ continuous in $r=0$ such that
:
\begin{equation}
f(x)\simeq f(r\theta)=r^k (f_k(\theta)+R(r,\theta)).
\end{equation}
If $\theta_0\notin C(f_k)$ and $\theta$ is close to $\theta_0$  we
simply choose :
\begin{gather*}
(y_{2},...,y_{n})(r,\theta) =(\theta _{1},...,\theta _{n-1}), \\
y_{1}(r,\theta)=r|f_{k}(\theta )+R(r,\theta)|^{\frac{1}{k}}.
\end{gather*}
In these coordinates the phase becomes $y_{1}^{k}$ if
$f_{k}(\theta_0)$ is positive (resp. $-y_{1}^{k}$ for a negative
value) and the Jacobian satisfies $|J y|(0,\theta )=|f_{k}(\theta
)|^{\frac{1}{k}}\neq 0$ locally. Now, let $\theta_0\in C(f_k)$. Up to a permutation, we can suppose that
$\partial_{\theta_{1}}f_{k}(\theta_{0})\neq 0$. We accordingly choose the
new local coordinates :
\begin{gather*}
(y_{1},y_{3},...,y_{n})(r,\theta)=(r,\theta
_{2},...,\theta _{n-1}), \\
y_{2}(r,\theta)=f_{k}(\theta)+R(r,\theta).
\end{gather*}
Since we have
$|Jy|(0,\theta_0)=|\partial_{\theta_{1}}f_{k}(\theta_{0})|\neq 0$,
lemma follows. \hfill{$\blacksquare$}
\begin{remark}\rm{These coordinates are
admissible and this leads to an adapted system of charts via a
partition of unity on $\mathbb{S}^{n-1}$ introduced below. Also
the condition $(H_1)$ insures the existence of a canonical measure on $C(f_k)$.}
\end{remark}
To use Lemma \ref{normal form} we introduce an adapted partition of
unity on $\mathbb{S}^{n-1}$. We pick cut-off functions
$\Psi_{j}\in C_0^{\infty}(\mathbb{S}^{n-1})$, $0\leq
\Psi_j \leq 1$, $\sum \Psi_j =1$ in a tubular neighborhood of $C(f_k)$,
with supports chosen so that normal forms of Lemma \ref{normal
form} exist, for $r$ small enough, in a conic neighborhood of
$\mathrm{supp}(\Psi_j)$. By compactness this set of functions can
be chosen finite and we obtain a partition of unity on
$\mathbb{S}^{n-1}$ by adding $\Psi_{0}=1-\sum \Psi_{j}$ to our
family. The support of $\Psi_{0}$ is not connected and we
define $\Psi_{0}^{+}$, with $f_k(\theta)>0$ on
$\mathrm{supp}(\Psi_{0}^{+})$, and similarly we define
$\Psi_{0}^{-}$ where $f_k<0$, so that
$\Psi_{0}=\Psi_{0}^{+}+\Psi_{0}^{-}$. If we accordingly split up
$I_s(z)$ we obtain :
\begin{gather*}
I_{s}^{\pm}(z)=\int\limits_{\mathbb{R\times R}_{+}\times
\mathbb{S}^{2n-1}} \Psi_{0}^{\pm}(\theta) g( zf(r\theta),r\theta )\Theta(r\theta)r^{n-1}drd\theta\\
=\int\limits_{\mathbb{R}_{+}} G_{0}^{\pm}(\pm z y_1^k, y_{1})dy_{1},
\end{gather*}
respectively for the directions where $f_k(\theta)>0$ and
$f_k(\theta)<0$, also :
\begin{gather*}
I_{s}^{0,j}(z)=\int\limits_{\mathbb{R\times R}_{+}\times
\mathbb{S}^{n-1}} \Psi_{j}(\theta) g(z f(r\theta) ,r\theta )\Theta(r\theta)r^{n-1}drd\theta \\
=\int\limits_{\mathbb{R}_{+}\times \mathbb{R}}
G_{j}(zy_1^ky_2, y_1,y_2)dy_1 dy_2,
\end{gather*}
for the set $C(f_k)$. The new symbols are
respectively given by :
\begin{gather}
G_{0}^{\pm}(t,y_{1})=\int y^{\ast
}(\Psi_{0}^{\pm}(\theta)g(t,r\theta
)\Theta (r\theta)r^{n-1}|Jy |)dy_{2}...dy_{n}, \label{ampli1}\\
G_{j}(t,y_{1},y_{2})=\int y ^{\ast }(\Psi_{j}(\theta)g(t,r\theta
)\Theta (r\theta)r^{n-1}|Jy |)d y_{3}...dy_{n}. \label{ampli2}
\end{gather}
Hence, the singular part of our integral can be written as a
finite sum :
\begin{equation}
I_s(z)=I_s^{-}(z)+I_s^{+}(z)+\sum\limits_{j} I_s^{0,j}(z),
\end{equation}
where each term of the r.h.s. will be treated by elementary
methods. Note that $I_s^{-}(z)$ and $I_s^{+}(z)$ can be treated as in the previous section.
\begin{remark}\label{degres amplitude}
\rm{Since $y_{1}(r,\theta ) =r$, our new symbols satisfy
$G_{j}(t,y_1,y_2)=\mathcal{O}(y_1^{n-1})$, near $y_1=0$. Since
the asymptotic expansion involves delta-Dirac distributions w.r.t $y_1$, cf. Lemma \ref{Theo IO 2eme carte} below,
the dimension will cause a shift in the expansion.}
\end{remark}
For $a\in C_0^{\infty}(\mathbb{R}\times \mathbb{R}_+\times\mathbb{R} )$, we define the family of elementary fiber integrals :
\begin{gather}
I_{n,k}(z)=\int\limits_{0}^{\infty }(\int\limits_{\mathbb{R}}a(zy_1^k y_2,y_1,y_2)dy_2)y_1^{n-1}dy_1.
\end{gather}
\begin{lemma}\label{Theo IO 2eme carte}
There exists a sequence of distributions $(D_{j,p})$ such that :
\begin{equation}
I_{n,k}(z)\sim \sum\limits_{p=0,1} \sum\limits_{j\in\mathbb{N},\text{ } j\geq n} D_{j,p}(a) z^{-\frac{j}{k}} \log(z)^p,
\text{ as } z\rightarrow \infty,
\end{equation}
where the logarithms only occur when $(j/k)$ is an integer. As concerns the leading term,
if $(n/k)\notin \mathbb{N}^{*}$ we obtain :
\begin{equation}\label{equation for I(nk)}
I_{n,k}(z)= z^{-\frac{n}{k}}d(a)+\mathcal{O}(z
^{-\frac{n+1}{k}}\mathrm{log}(z)),
\end{equation}
with :
\begin{gather*}
d(a)=C_{n,k} \int\limits_{0}^{\infty}\int\limits_{0}^{\infty} t^{\frac{n}{k}-1}
y_2^{n-\frac{n}{k}} \left( \partial^{n}_{y_2} a(t,0,y_2)+\partial^{n}_{y_2} a(-t,0,-y_2)\right) dy_2dt.
\end{gather*}
But when $n/k=p\in\mathbb{N}^{*}$, we have :
\begin{equation*}
I_{n,k}(z)= \frac{1}{k} z^{-\frac{n}{k}}\log(z)
\int\limits_{\mathbb{R}} |t|^{p-1} \partial_{y_2}^{p-1} a(t,0,0)dt
+\mathcal{O}(z^{-\frac{n}{k}}).
\end{equation*}
\end{lemma}
\begin{remark}\rm{The remainder of Eq. (\ref{equation for I(nk)}) can be optimized to $\mathcal{O}(z
^{-\frac{n+1}{k}})$ when $(n+1)/k$ is not an integer, as shows the proof below.}
\end{remark}
\noindent\textit{Proof.} By a standard density argument we can assume that the amplitude is of the form
$a(s,y_1,y_2)=f(s)b(y_1,y_2)$. The justification is that our coefficients below are computed by
continuous linear functionals, i.e. distributions. We define the Melin transforms of $f$ as :
\begin{equation}
M_{\pm}(\xi)= \int\limits_{0}^{\infty} s^{\xi-1} f(\pm s)ds.
\end{equation}
We split-up $I_{n,k}$ as $I_{+}$ and $I_{-}$ by separating integrations $y_2>0$ and $y_2<0$. Via Melin's
inversion formula, we accordingly obtain :
\begin{equation}\label{analytic1}
I_{+}(z)=
\frac{1}{2i\pi }\int\limits_{\gamma} M_{+}(\xi) z^{-\xi}
\int\limits_{\mathbb{R}_{+}^2} (y_1 y_2^{k})^{-\xi} b(y_1,y_2) y_1^{n-1} dy_1dy_2 d\xi,
\end{equation}
where $\gamma=c+i\mathbb{R}$ and $0<c<k^{-1}$. Similarly we have :
\begin{equation}\label{analytic2}
I_{-}(z)=
\frac{1}{2i\pi }\int\limits_{\gamma} M_{-}(\xi) z^{-\xi}
\int\limits_{\mathbb{R}_{+}^2} (y_1 y_2^{k})^{-\xi} b(y_1,-y_2) y_1^{n-1} dy_1dy_2 d\xi,
\end{equation}
The existence of a full asymptotic expansion is a direct consequence of :
\begin{lemma}\label{poles extensions}
The family of distributions $\xi\mapsto (y_1 y_2^{k})^{-\xi}$ on $C_0^{\infty}(\mathbb{R}_{+}^2)$ initially defined
in the domain $\Re(\xi)< k^{-1}$ is meromorphic on $\mathbb{C}$ with poles :
$\xi_{j,k} = j/k$, $j\in\mathbb{N}^{*}$. These poles are of order 2 when $\xi_{j,k}\in \mathbb{N}^*$ and of order 1
otherwise.
\end{lemma}
\textit{Proof.} We form the Bernstein-Sato polynomial $b_k$ attached to our problem :
\begin{gather*}
T(y_2y_1^{k})^{1-\xi}:=\frac{\partial}{\partial y_2}
\frac{\partial^{k}}{\partial y_1^{k}}
(y_2y_1^{k})^{1-\xi}=b_k(\xi)(ty_2y_1^{k})^{-\xi}, \\
b_k(\xi)=(1-\xi)\prod\limits_{j=1}^{k}(j-k\xi).
\end{gather*}
If $\Re(\xi)<k^{-1}$, $(k+1)$-integrations by parts yield :
\begin{equation*}
\int\limits_{\mathbb{R}_{+}^2} (y_1 y_2^{k})^{-\xi} f(y_1,y_2)dy_1dy_2=\frac{(-1)^{k+1}}{b_k(\xi)}
\int\limits_{\mathbb{R}_{+}^2} (y_1 y_2^{k})^{1-\xi} (Tf)(y_1,y_2)dy_1dy_2.
\end{equation*}
Now the integral in the r.h.s. is analytic in $\Re (\xi)<1+k^{-1}$. After $m$ iterations
the poles, with their orders, can be read off the rational functions :
\begin{equation}
\mathfrak{R}_m(\xi) =\prod\limits_{p=1}^m \frac{1}{b_k(\xi-p)}.
\end{equation}
This gives the result since $m$ can be chosen arbitrary large. $\hfill{\blacksquare}$
\begin{remark}\rm{Starting from Eqs. (\ref{analytic1},\ref{analytic2}) it is possible
to obtain an asymptotic expansion as $z\rightarrow 0^{+}$.
By shifting the path $\gamma$ to the left, only the poles of the Melin transforms
contribute and we obtain an expansion in power of $z$.
This can be predicted by a Taylor expansion of the associated oscillatory integral but the improvement
is that one can establish sharp integral remainders.}
\end{remark}
Hence, the following functions are meromorphic on $\mathbb{C}$ :
\begin{equation}
\mathfrak{g}^{\pm}(\xi)=\int\limits_{\mathbb{R}_{+}^2} (y_1 y_2^{k})^{-\xi} b(y_1,\pm y_2)dy_1dy_2.
\end{equation}
A classical result, see e.g. \cite{Bl-Hand}, is that $M_{\pm}(c+ix)\in \mathcal{S}(\mathbb{R}_x)$
when $c \notin -\mathbb{N}$. If we shift the path of integration $\gamma$ to the right in our
integral representation, Cauchy's residue method provides the asymptotic expansion.
In fact for any $d>c$, outside of the poles, we have :
\begin{gather*}
I_{+}(z)=\int\limits_{c+i\mathbb{R}} z^{-\xi} M_{+}(\xi) \mathfrak{g}^{+}(\xi) d\xi\\
= \sum\limits_{c<\xi_{j,k}<d} \mathrm{res}(z^{-\xi}M_{+}\mathfrak{g})(\xi_{j,k})+
\int\limits_{d+i\mathbb{R}} z^{-\xi} M_{+}(\xi) \mathfrak{g}^{+}(\xi) d\xi.
\end{gather*}
Since $d$ is not a pole the last integral can be estimated via :
\begin{equation}
|\int\limits_{d+i\mathbb{R}} z^{-\xi} M_{+}(\xi) \mathfrak{g}^{+}(\xi) d\xi|\leq C(f,b) z^{-d}=\mathcal{O}(z^{-d}),
\end{equation}
where, for each $d$, the constant $C$ involves the $L^1$-norm of a finite number derivatives of $b$.
This will indeed lead to an asymptotic expansion with precise remainders.
Applying this method to $I_{+}(z)$ and $I_{-}(z)$ we obtain the existence of a full asymptotic expansion of the form :
\begin{equation}
I(z) \sim \sum\limits_{p=0,1} \sum\limits_{j\in\mathbb{N}^*} C_{j,p} z^{-\frac{j}{k}} \log(z)^p.
\end{equation}
Moreover, by Lemma \ref{poles extensions}, these logarithms only occur when $j/k$ is integer.\medskip\\
Now, we compute the leading term of this expansion for our particular problem. To avoid unnecessary discussions and
calculations below, we remark that we can commute the polynomial weight via :
\begin{gather}
T ((y_1 y_2)^{1-\xi} y_1^{n-1})=\mathfrak{b}(\xi)(y_2y_1^{k})^{-\xi}y_1^{n-1}, \\
\mathfrak{b}(\xi)=(1-\xi)\prod\limits_{j=1}^{k}(j-k\xi+n-1).
\end{gather}
By iteration, we obtain that the poles are the rational numbers :
\begin{equation*}
\xi_{p,j,k,n}=p+\frac{j+n-1}{k},\text{ } j\in[1,...,k], \text{ } p\in\mathbb{N}.
\end{equation*}
For all $\alpha-\beta>1$, $\beta\in\mathbb{N}^{*}$, we have :
\begin{equation*}
\int\limits_{0}^{\infty} \partial_r^{\beta}  (r^\alpha f(r)dr) =0, \text{ }\forall f\in C_{0}^{\infty}.
\end{equation*}
If we apply this to the integral w.r.t. $y_1$, there is no contribution before :
\begin{equation}
\xi_{0}= \frac{n}{k}.
\end{equation}
To compute the first effective residue we must distinguish out
the case where $\xi_0$ is an integer or not. The optimal number of iterations to reach $\xi_0$
is $\mathrm{E}(n/k)+1$ but, by analytic continuation, any number bigger than this one is acceptable.
A fortiori we can use $n$ iterations and our starting point will be :
\begin{equation}
z^{-\xi}M_{+}(\xi) (-1)^{n(k+1)} \mathfrak{B}_n (\xi)
\int\limits_{\mathbb{R}_{+}^{2}} (y_1^k y_2)^{n-\xi} y_1^{n-1} T^n b(y_1,y_2) dy_1dy_2,
\end{equation}
with :
\begin{equation}
\mathfrak{B}_n(\xi)=\prod\limits_{l=0}^{n-1} \frac{1}{\mathfrak{b}(\xi-l)}.
\end{equation}
\textbf{Case of $\xi_0$ simple pole.}\\
In this case our residue is simply given by :
\begin{equation}
C z^{-\frac{n}{k}} M_{+}(\frac{n}{k})
\int\limits_{\mathbb{R}_{+}^{2}} (y_1^k y_2)^{n-\frac{n}{k}} y_1^{n-1} T^n b(y_1,y_2) dy_1dy_2,
\end{equation}
with :
\begin{equation}
C =\lim\limits_{\xi\rightarrow \frac{n}{k}} (-1)^{n(k+1)}(\xi-\frac{n}{k})\mathfrak{B}_n(\xi).
\end{equation}
In particular we can compute the integral w.r.t. $y_1$ via :
\begin{equation*}
\int\limits_{0}^{\infty} y_1^{kn-1} \partial^{kn}_{y_1}  (\partial^{n}_{y_2} b(y_1,y_2)) dy_1
=(-1)^{kn}(kn-1)!\partial^{n}_{y_2} b(0,y_2).
\end{equation*}
A similar result holds for $I_{-}$ and we obtain :
\begin{gather}
I_{+}(z)=z^{-\frac{n}{k}} C_{n,k}  M_{+}(\frac{n}{k})
\int\limits_{0}^{\infty} y_2^{n-\frac{n}{k}} (\partial^{n}_{y_2} b)(0,y_2)dy_2+R_1(z),\\
I_{-}(z)=z^{-\frac{n}{k}} C_{n,k}  M_{-}(\frac{n}{k})
\int\limits_{0}^{\infty} y_2^{n-\frac{n}{k}} (\partial^{n}_{y_2} b)(0,-y_2)dy_2+R_2(z).
\end{gather}
Here $C_{n,k}$ is the canonical constant :
\begin{equation}\label{canonical regularization}
C_{n,k}=\frac{1}{k} \prod\limits_{j=1}^n \frac{-1}{j-\frac{n}{k}}.
\end{equation}
Also, according to the analysis above, each remainder is of order $\mathcal{O}(z^{-\frac{n+1}{k}})$
if $(n+1)/k \notin \mathbb{N}$ and $\mathcal{O}(z^{-\frac{n+1}{k}}\log(z))$
otherwise.\medskip\\
\textbf{Case of $\xi_0$ double pole.}\\
If $h$ is meromorphic with a pole of order 2 in $\xi_0$ we have :
\begin{equation*}
\mathrm{res} (h)(\xi_0)=\frac{1}{2} \lim\limits_{\xi\rightarrow \xi_0} \frac{\partial}{\partial \xi}
(\xi-\xi_0)^2 h(\xi).
\end{equation*}
Applying this principle to our residue we obtain, via Leibnitz's rule, that :
\begin{equation}
I_{+}(z)=B \log(z)z^{-\frac{n}{k}} +\mathcal{O}(z^{-\frac{n}{k}}).
\end{equation}
We can compute the distribution $B$ as before and we find :
\begin{gather*}
B=-\frac{1}{2} D_{n,k} M_{+}(\frac{n}{k})
\int\limits_{0}^{\infty} y_2^{n-\frac{n}{k}} (\partial^{n}_{y_2} b)(0,y_2)dy_2,\\
D_{n,k}=(-1)^{n(k+1)} \lim\limits_{\xi \rightarrow \frac{n}{k}} (\xi-\frac{n}{k})^2 \mathfrak{B}_n(\xi).
\end{gather*}
Since $p=n/k$ is an integer, by integration by parts we obtain :
\begin{equation}
\int\limits_{0}^{\infty} y_2^{n-p} (\partial^{n}_{y_2} b)(0,y_2)dy_2
=(-1)^{n-p+1}(n-p)!  \partial_{y_2}^{p-1} b(0,0).
\end{equation}
Since a similar result holds for $I_{-}$, we obtain the desired result by gathering all the constants and summation.
Finally, we can extend our formulas since all coefficients in the expansion are of the form :
\begin{equation*}
\left\langle T^j ,f\otimes b\right\rangle=\left\langle T^j_1,f\right\rangle \left\langle T^j_2,b\right\rangle,
\text{ }T^j_{1,2}\in \mathcal{D}'.
\end{equation*}
By linearity and continuity, the result holds for a symbol $a(t,y_1,y_2)$.
\hfill{$\blacksquare$}\medskip\\
Taking Remark \ref{degres amplitude} into account, to avoid
unnecessary calculations we define :
\begin{gather}
G_{0}^{\pm}(t,y_1)=y_1^{n-1}\tilde{G}_{0}^{\pm}(t,y_1),\\
G_j(t,y_1,y_2)=y_1^{n-1}\tilde{G}_{j}(t,y_1,y_2).
\end{gather}
\textbf{Directions where $f_k(\theta)\neq 0$.}\\
By Lemma \ref{lemma for extremum}, the first non-zero coefficient, obtained for $l=n-1$, is :
\begin{equation*}
\frac{z^{-\frac{n}{k}}}{k}\frac{1}{(n-1)!}\left\langle
t_{+}^{\frac{n-k}{k}}\otimes \delta
_{0}^{(n-1)},G_{0}^{+}(t,y_{1})\right\rangle =\frac{z^{-\frac{n}{k}}}{k}\int\limits_{\mathbb{R}}
t_{+}^{\frac{n-k}{k}}
\tilde{G}_{0}^{+}(t,0)dt.
\end{equation*}
By construction, we have :
\begin{equation*}
\tilde{G}_{0}^{+}(t,0)=\int\limits_{\mathbb{S}^{n-1}}g(t,0)\Psi_{0}^{+}(\theta)|f_{k}(\theta
)|^{-\frac{n}{k}}d\theta.
\end{equation*}
A similar computation gives the contribution of
$\mathrm{supp}(\Psi_{0}^{-})$, and we obtain :
\begin{gather} \label{positive}
I_s^{+}(z)= z^{-\frac{n}{k}} \left\langle
t_{+}^{\frac{n-k}{k}}, g(t,0) \right\rangle\frac{1}{k}
\int\limits_{\mathbb{S}^{n-1}}\Psi_{0}^{+}(\theta)|f_{k}(\theta
)|^{-\frac{n}{k}}d\theta+\mathcal{O} (z^{-\frac{n+1}{k}}),\\
\label{negative}%
I_s^{-}(z) = z^{-\frac{n}{k}} \left\langle
t_{-}^{\frac{n-k}{k}}, g(t,0) \right\rangle
\frac{1}{k}\int\limits_{\mathbb{S}^{n-1}}\Psi_{0}^{-}(\theta)|f_{k}(\theta
)|^{-\frac{n}{k}}d\theta +\mathcal{O} (z^{-\frac{n+1}{k}}).
\end{gather}
\textbf{Microlocal contribution of $C(f_k)$.}\newline%
Here we examine the contribution of terms $I_{s}^{0,j}(z)$. According to the analysis above,
we must distinguish out the case $k$ divides $n$.\medskip\\
\textit{a) Case of $k>n$, non-integrable singularity on $\mathfrak{S}$.}\\
Here $n/k \in ]0,1[$, so that the singularity on the blow-up is integrable. Via Lemma \ref{Theo IO 2eme carte}, the contribution of $I_{s}^{0,j}(z)$
is given by :
\begin{equation*}
\frac{1}{k}z^{-\frac{n}{k}}\int\limits_{\mathbb{R}_{+}^2} |t|^{\frac{n}{k}-1}
|y_2|^{-\frac{n}{k}} \left(\tilde{G}_{j}(t,0,y_2)+\tilde{G}_{j}(-t,0,-y_2)\right)dt dy_2
+\mathcal{O}(z^{-\frac{n+1}{k}}\mathrm{log}(z)).
\end{equation*}
Reminding that $y_2(t,0,\theta)=f_k(\theta)$, we obtain :
\begin{equation*}
\int\limits_{\mathbb{R}_{+}}
|y_2|^{-\frac{n}{k}} \tilde{G}_{j}(t,0,y_2)dy_2=g(t,0)\int\limits_{\{f_k(\theta)\geq0\}}|f_k
(\theta) |^{-\frac{n}{k}}\Psi_j(\theta) d\theta.
\end{equation*}
Since $\Psi_0^{\pm}$, $\Psi_j$ is a partition of unity on $\mathbb{S}^{n-1}$, by summation of all local
contributions $I_s(z)$ is asymptotically equivalent to :
\begin{equation*}
\frac{z^{-\frac{n}{k}}}{k}  \left ( \left\langle
t_{+}^{\frac{n}{k}-1},g(t,0) \right\rangle
\int\limits_{\{ f_k \geq 0\}} |f_k
(\theta) |^{-\frac{n}{k}}d\theta
+ \left\langle t_{-}^{\frac{n}{k}-1},g(t,0)\right\rangle
\int\limits_{\{ f_k \leq 0\}} |f_k (\theta)
|^{-\frac{n}{k}}d\theta \right ).
\end{equation*}
Note that none of these coefficients are
equal unless $g$ or $f_k$ are symmetric. $\hfill{\blacksquare}$\medskip\\
\textit{b) Case of $p=n/k$ integer.}\\
Here the contribution of $I_{s}^{0,j}(z)$ is dominant and we obtain :
\begin{equation*}
I_{s}^{0,j}(z)\sim \frac{1}{k}\log (z)z^{-p}
\int\limits_{\mathbb{R}} |t|^{p-1} \partial^{p-1}_{y_2}
\tilde{G}_j(t,0,0)dt+\mathcal{O}(z^{-p}).
\end{equation*}
Unless $p=1$, there is no way to take the limit directly, and the geometric properties are still hidden in the Jacobian.
But we will reach the result by the Schwartz kernel technic. Clearly, it is enough to evaluate our derivative
and to integrate w.r.t. $t$. With $s=(s_1,s_2)\in\mathbb{R}^2$, we write the evaluation as :
\begin{equation*}
\partial^{p-1}_{y_2} \tilde{G}_j(t,0,0)=\frac{1}{(2\pi)^2} \int e^{i\left\langle
s,(y_1,y_2)\right\rangle} (is_2)^{p-1}
\tilde{G}_j(t,y_1,y_2) dy_1dy_2 ds.
\end{equation*}
Here we have used an oscillatory Schwartz kernel for $\delta_{y_1} \otimes \delta^{p-1}_{y_2}$.
This integral representation allows to inverse our diffeomorphism to obtain :
\begin{equation*}
\partial^{p-1}_{y_2}
\tilde{G}_j(t,0,0)=\frac{1}{(2\pi)^2} \int e^{i\left\langle
s,(r,y_2(r,\theta))\right\rangle} (i s_2)^{p-1} t^{p-1}
g(t,r\theta) \Psi_j(\theta) dtdr d\theta ds.
\end{equation*}
Extending the integrand by 0 for $r<0$, the normalized integral w.r.t. $(r,s_1)$ provides $\delta_r$.
By construction $y_2(0,\theta)=f_k(\theta)$, hence :
\begin{equation}
\partial^{p-1}_{y_2}
\tilde{G}_j(t,0,0)=g(t,0)\frac{1}{(2\pi)} \int\limits_{\mathbb{R}\times\mathbb{S}^{n-1}}
e^{i u f_k(\theta)} (i u)^{p-1} \Psi_j(\theta) d\theta du .
\end{equation}
This Fourier integral makes sense with $\mathbb{S}^{n-1}$ compact. We recall the density :
\begin{equation}
J_j(w)=\int \limits_{\{f_k(\theta)=w\}} \Psi_j(\theta) dL_w(\theta),
\end{equation}
where $dL_w$ is the density induced by the Leray-form $dL_{f_k}$ : $d f_k \wedge dL_{f_k}(\theta)=d\theta$.
Note that all these objects can be constructed by mean of local coordinates under the only condition
that $\mathrm{supp}(\Psi_j)$ is small enough near $C(f_k)$. Since $f_k$ is continuous on $\mathbb{S}^{n-1}$,
each $J_j(w)$ is compactly supported and smooth near the origin. The sum
over all the $\Psi_j$ gives the geometric contribution :
\begin{equation}
\frac{1}{(2\pi)} \int\limits_{\mathbb{R}^2} e^{i uw} (i u)^{p-1} \sum\limits_j J_j(w) dwdu=\frac{d^{p-1}\mathrm{Lvol}}{d w^{p-1}}(0).
\end{equation}
By integration w.r.t. $t$ we obtain the result for a double pole in general position.
\begin{remark}
\rm{The case $n=k$ is directly accessible. In this case we obtain :
\begin{equation*}
I_s(z)=\frac{1}{k} \frac{\log(z)}{z} \mathrm{LVol}(0)\int\limits_{\mathbb{R}} g(t,0)dt+\mathcal{O}(z^{-1}).
\end{equation*}
Here $\mathrm{LVol}(0)$ is the Liouville volume of $C(f_k)$. Note that $I_s$ dominates $I_r$.}
\end{remark}
\textit{c) $k<n$ and simple pole, integrable singularity on $\mathfrak{S}$.}\\
Finally, we treat the case of a simple pole with a non-integrable singularity on $\mathbb{S}^{n-1}$.
For the positive part, the distributional coefficients are :
\begin{equation}
\left \langle \mu^j_{+}, g\right\rangle = C_{n,k} \int\limits_{0}^{\infty} \int\limits_{0}^{\infty}
t^{\frac{n}{k}-1} y_2^{n-\frac{n}{k}} (\partial^{n}_{y_2} \tilde{G}_j)(t,0,y_2)dy_2dt.
\end{equation}
Clearly, we can use the same globalization technic as above.
The sum of all $\mu^j_{+}$, completed with the main term of $J_s^{+}$, provides :
\begin{equation}
\left\langle T_{+}, g\right\rangle =C_{n,k}\int\limits_{0}^{\infty}|t|^{\frac{n}{k}-1}g(t,0)dt
\int\limits_{0}^{\infty} w^{n-\frac{n}{k}} \frac{ \partial^{n} L(w)}{\partial w^n}  dw,
\end{equation}
where the integral w.r.t. $w$ is absolutely convergent since the measure $L(w)$ is compactly supported.
A similar result holds for the directions where $f_k(\theta)<0$
and we obtain the result stated in Theorem \ref{main2}.$\hfill{\blacksquare}$\medskip\\
Now we detail the construction of the distributional bracket of the part c) of Theorem \ref{main2}.
Let be $\chi\in C_0^{\infty}$, $0\leq \chi \leq 1$ on $\mathbb{R}$, chosen such that
$\chi=1$ near the origin and $\chi(u)=0$ for $|u|\geq \varepsilon$. If $\varepsilon>0$ is small enough,
$\mathrm{Lvol}$  is smooth on $]-\varepsilon, \varepsilon[$.
We write the geometric contribution as :
\begin{equation*}
\left\langle T,\mathrm{Lvol} \right\rangle= \left\langle T, \chi \mathrm{Lvol} \right\rangle
+\left\langle T, (1-\chi)\mathrm{Lvol}  \right\rangle
\end{equation*}
Away from the origin, e.g. for $u>0$, we obtain directly :
\begin{equation}
C_{n,k}\left\langle \frac{d^n}{du^n} u_{+}^{n-\frac{n}{k}}, (1-\chi)(u)\mathrm{Lvol}(u) \right\rangle
=\frac{1}{k}\int\limits_{\{f_k(\theta) >0\}} (1-\chi)(f_k(\theta))|f_k(\theta)|^{-\frac{n}{k}}d\theta.
\end{equation}
Note that Eq. (\ref{canonical regularization}) for $C_{n,k}$ justifies the
normalization of Eq. (\ref{normalized derivatives}).
On $\mathrm{supp}(\chi)$, we use the local regularity of $\mathrm{Lvol}(u)$ near $u=0$ and
integrations by parts to conclude. Finally, to obtain a totally rigorous treatment, we remark that
since $C(f_k)$ is compact and $f_k$ is continuous we can choose our partition of unity such that $\sum\Psi_j=1$ for
$|f_k(\theta)|\leq 2\varepsilon$. $\hfill{\blacksquare}$
\begin{remark}\rm{The key point here is that we can put in duality the
distributions $\partial^n_u |u|_{\pm}^{n-\alpha}$ and $\mathrm{Lvol}(u)$ since their singular supports are disjoints.}
\end{remark}
\section{Elementary examples.}
\textbf{A family of extremum.} If $p$ is an even integer and $a=(a_1,...a_n)\in (\mathbb{R}_{+}^{*})^n$,
we consider a $n$ dimensional integral obtained by integration in the fiber of :
\begin{equation}
f_{a,p}(x)=\sum \limits_{j=1}^n a_j x_j^p
\end{equation}
For example :
\begin{equation*}
\int\limits_{\mathbb{R}^n} e^{-z f_{a,p}(x)}dx=
\prod\limits_{j=1}^{n} (\int\limits_{\mathbb{R}} e^{-z a_j x^p} dx)=(2 \Gamma
(1+\frac{1}{p}))^n z^{-\frac{n}{p}}\prod\limits_{j=1}^{n}a_j^{-\frac{1}{p}},
\end{equation*}
with $g(t,x)=e^{-t}$ symbol of order 0 on $\mathbb{R}^n$. Since we have here a local minimum, the
distributional factor is simply :
\begin{equation}
\int\limits_{0}^{\infty} t^{\frac{n-p}{p}} e^{-t}dt= \Gamma(n/p).
\end{equation}
The Jacobian of the standard polar coordinates gives the value of the generalized elliptic
integrals :
\begin{equation} \label{elliptic}
\int\limits_{\mathbb{S}^{n-1}}
(f_{a,p}(\theta))^{-\frac{n}{p}}d\theta=\frac{1}{p}\frac{(2 \Gamma
(1+\frac{1}{p}))^n}{\Gamma(n/p) } \prod\limits_{j=1}^{n}a_j^{-\frac{1}{p}}.
\end{equation}
Note that Eq. (\ref{elliptic}) can be analytically continued with a phase $\sum a_j|x_j|^\alpha$, $\alpha>0$.\medskip\\
\textbf{Conical singularities.} In 2 dimensions
let be :
\begin{equation*}g(t,u)=\frac{e^{-||u||^2}}{(1+t^2)}.
\end{equation*}
With $f(x,y)=x^2-y^2$ we have $\mathfrak{S}=\{(x,x)\} \cup \{(-x,x)\}$. $\mathfrak{S}$ is not compact but
the exponential decrease of $g$ will compensate. Accordingly let be :
\begin{equation}
I(z)=\int\limits_{\mathbb{R}^2} \frac{e^{-(x^2+y^2)}}{1+z^2
(x^2-y^2)^2}dxdy.
\end{equation}
Since $n/k=1$, we expect a logarithm for the leading term. Since :
\begin{equation*}
\int\limits_0^{2\pi} \frac{d\theta}{1+\alpha \cos^2
(2\theta)}=\frac{2\pi}{\sqrt{1+\alpha}},\text{ } \forall \alpha
>-1,
\end{equation*}
by passage in polar coordinates we obtain :
\begin{equation*}
I(z)=2\pi \int\limits_{0}^{\infty} e^{-r^2} \frac{rdr}{\sqrt{1+z^2
r^4}}=\pi \int\limits_{0}^{\infty} e^{-u} \frac{du}{\sqrt{1+z^2
u^2}} \sim \pi \frac{\log(z)}{z}.
\end{equation*}
Always in dimension 2, a singularity of degree 4 will generally not give a logarithmic leading term.
For example, we have :
\begin{equation*}
\tilde{I}(z)=\int\limits_{\mathbb{R}^2}\frac{e^{-(x^2+y^2)}}{1+z^2(x^4-y^4)^2}dxdy
=\pi\int\limits_{0}^{\infty} e^{-u} \frac{du}{\sqrt{1+z^2
u^4}} \sim \frac{\pi}{\sqrt z} \int\limits_{0}^{\infty} \frac{du}{\sqrt{1+u^4}}.
\end{equation*}

\end{document}